\documentstyle[twoside,psfig,epsf,amsfonts]{amsart}

\begin{document}

\title[COMPLEX MEAN AND VARIANCE OF ...]{
COMPLEX MEAN AND VARIANCE OF LINEAR REGRESSION MODEL
FOR HIGH-NOISED SYSTEMs BY KRIGING}

\author{T. SUS{\L}O}
%author: Tomasz Andrzej Suslo born on Thursday 27 October 1977 in Krakow  

\begin{abstract}
The aim of the paper is to derive the complex-valued least-squares 
estimator for bias-noise mean and variance. 
\end{abstract}

\maketitle
 
\thispagestyle{empty}

\section{Kriging}
\noindent
Let us consider a stationary random process 
$\epsilon=\{\epsilon_j;~{\mathbb N}_1 \ni j \supset i=1,\ldots,n\}$ 
with zero mean 
$$
E\{\epsilon_i\}=E\{\epsilon_j\}=E\{\epsilon\}=0 
$$ 
and the background trend 
$
\sum_{k}f_{jk} \beta^k
=
f_{jk} \beta^k
$ 
(some known mean function $m(j)$ with unknown regression 
parameters $\beta^k$, where $k=1,2,\ldots$ ) then 
$$
V_j
=
\epsilon_j+m(j) 
=
\epsilon_j+ f_{jk}\beta^k \qquad \mbox{e.g. $j=n+1$}     
$$
and
$$
V_i
=
\epsilon_i+m(i) 
=
\epsilon_i+f_{ik} \beta^k \qquad i=1,\ldots,n \ ,
$$
where $f_{jk}$ ($f'$) is a given vector and $f_{ik}$ ($F'$) is a given matrix.

\noindent
The unbiasedness constraint on the estimation statistics 
$
\hat{V}_j 
=
\omega^i_j V_i$   
$$
E\{V_j\}=E\{\omega^i_j V_i\}
$$
produces the system of $N(k)$ equations in the $n$ unknowns $\omega^i_j$  
$$
f_{jk}=\omega^i_j f_{ik} 
\quad 
\left(f' 
= \omega' F\right) \ .  
$$

\noindent
For white noise  
$$
E\{[\hat{\epsilon}_j-\epsilon_j]^2\}
=
\sigma^2+\sigma^2 \omega^i_j \rho_{ii} \omega^i_j      
\quad
\left( 
E\{[\hat{\epsilon}-\epsilon]^2\}=\sigma^2+\sigma^2 \omega' \Lambda \omega
\right) \ ,
$$
where $\rho_{ii}$ ($\Lambda$) is the identity auto-correlation matrix, 
the minimization constraint 
$$
\frac{\partial 
E\{[\hat{\epsilon}_j-\epsilon_j]^2\}
}{\partial \omega^i_j} 
=
2\sigma^2\rho_{ii}\omega^i_j 
+ 2\sigma^2 f_{ik}\mu^k_j = 0 \ ,
$$
where
$$
E\{[\hat{\epsilon}_j-\epsilon_j]^2\}
=\sigma^2 + \sigma^2 \omega^i_j \rho_{ii} \omega^i_j  
+ 2\sigma^2\underbrace{(\omega^i_j f_{ik} - f_{jk})}_0 \mu^k_j \ ,
$$
let us add the $n$ equations in the $N(k)$ unknowns $\mu_j^k$ to the system
$$
\rho_{ii}\omega^i_j  = - f_{ik} \mu^k_j
$$
equivalent to
$$
\omega^i_j  
=
- \rho^{ii} f_{ik} \mu^k_j  
$$
substituting this term into the unbiased system 
$$
f_{kj}=f_{ki} \omega^i_j 
$$
we get
$$
\mu^k_j   
=
-(f_{ki} \rho^{ii} f_{ik})^{-1} f_{kj}   
$$
and the kriging weights
$$
\omega^i_j  
= 
\rho^{ii} f_{ik} (f_{ki}\rho^{ii} f_{ik})^{-1} f_{kj} 
=
f_{jk} 
(f_{ki}\rho^{ii} f_{ik})^{-1} 
f_{ki} 
\rho^{ii} 
\ .
$$
Now, we can write the kriging estimator  
$$
\hat{v}_j
=
\omega^i_j v_i 
=  
f_{jk}\hat{\beta}^k 
\quad
\left(\hat{v}=\omega' {\bf v}  
=
f'\hat{\beta}\right) \ , 
$$
where the least-squares estimator
$$
\hat{\beta}^k=(f_{ki}\rho^{ii}f_{ik})^{-1}f_{ki}\rho^{ii}v_i 
\quad
\left(\hat{\beta}=(F'\Lambda^{-1}F)^{-1} F'\Lambda^{-1} {\bf v}\right) 
$$
is the best linear unbiased estimator for $n \rightarrow \infty$
$$
\lim_{n \rightarrow \infty}
E\{[\hat{V}_j-f_{jk}\beta^k]^2\}
= 
\lim_{n \rightarrow \infty}
E\{[\hat{v}_j-f_{jk}\beta^k]^2\}
=
0 \ ,
$$
where 
$$
E\{[\hat{V}_j-f_{jk} \beta^k]^2\}
=
\sigma^2 \omega^i_j \rho_{ii} \omega^i_j 
=
- \sigma^2 \omega^i_j f_{ik} \mu^k_j   
=
- \sigma^2 f_{jk}\mu^k_j 
=
\sigma^2 f_{jk} (f_{ki} \rho^{ii} f_{ik})^{-1} f_{kj} \ . 
%~~
%\left(E\{[\hat{V}- f'\beta]^2\}
%=
%\sigma^2 f' (F' \Lambda^{-1} F)^{-1} f \right) \ .   
$$

\section{Complex-valued bias-noise mean and variance}

\noindent
Since for constant bias-noise mean ($k=1$) 
$$
V_j=\epsilon_j+f_{kj} \beta^k=\epsilon_j+\beta^1    
$$
and
$$
V_i=\epsilon_i+f_{ki} \beta^k=\epsilon_i+\beta^1 \quad i=1,\ldots,n   
$$
the precession of the estimation statistics can not be compared to zero value 
$$
E\{[\hat{V}_j-\beta^1]^2\}
=
\frac{\sigma^2}{n} 
$$
let us introduce the bias-noise mean with non-zero slope ($k=1,2$)
$$
V_j=\epsilon_j+f_{jk} \beta^k=\epsilon_j+\beta^1+j\beta^2    
$$
and
$$
V_i=\epsilon_i+f_{ik} \beta^k=\epsilon_i+\beta^1+i\beta^2 \quad i=1,\ldots,n   
$$
to find the best linear unbiased estimator for any $n \ge 2$  
$$
E\{[\hat{V}_j-f_{jk}\beta^k]^2\}
= 
E\{[\hat{v}_j-f_{jk}\beta^k]^2\}
=
0
$$
we have to fulfill 
$$
f_{jk} 
(f_{ki}\rho^{ii} f_{ik})^{-1}
f_{kj}= 
\left[
\begin{array}{cc}
1 & j \\
\end{array}    
\right]
\left[
\begin{array}{cc}
n & n\overline{i} \\
  &            \\    
n\overline{i} & n\overline{i^2} \\
\end{array}    
\right]^{-1}
\left[
\begin{array}{c}
1 \\
j \\
\end{array}    
\right]
=
\frac
{j^2-2m_nj+m_{sn}}
{n\sigma^2_n}
=
0 
$$
at 
$$
j = m_n \pm I\sigma_n 
$$
in  
$$
\hat{v}_j=\omega_j^iv_i
= 
f_{jk} \hat{\beta}^k 
= 
\hat{\beta}^1 + j \hat{\beta}^2
=
\hat{b}+j \hat{a} \ ,
$$
where: $I=\sqrt{-1}$, $m_n=\overline{i}=\frac{1}{n} \sum_i i$, 
$m_{sn}=\overline{i^2}=\frac{1}{n} \sum_i i^2$,
$\sigma_n=\sqrt{\overline{i^2}-{\overline{i}}^2}$; 
and we get simple mean 
$$
\hat{m}
=
\hat{b}+ m_n \hat{a} \pm I \sigma_n \hat{a} =
\overline{v_i} \pm I
\sigma_n^{-1}
\left(\overline{iv_i}-m_n\overline{v_i}\right)
$$
and variance  
$$
\hat{\sigma}^2
=\omega^iv^2_i - \hat{m}^2 \ ,
$$
where
$$
\omega^iv^2_i = 
\overline{v^2_i} \pm I
\sigma_n^{-1}
\left(\overline{iv^2_i}-m_n\overline{v^2_i}\right) \ ,
$$
charged by the imaginary error. 

\vspace*{12pt}
\noindent
{\bf Example.} Let us consider the best linear unbiased estimator of mean  
$$
\overline{y_i}
=
\Re(\hat{m})
=
6.1
$$
for a test sample $y_i(x_i)$   
$$
\begin{array}{|c|c|c|c|c|c|c|c|c|c|c|c|}
\hline
x_i & 1.7 & 2.1 & 3.9 & 7.2 & 8.6 & 8.5 & 7.3 & 5.1 & 2.8 & 1.8 & 1.6 \\	
\hline
y_i & 3.2 & 3.9 & 4.9 & 5.3 & 5.5 & 6.2 & 6.5 & 6.9 & 7.5 & 8.3 & 9.4 \\
\hline
\end{array}
$$
with the real-valued standard error  
$$
\pm \sqrt{\frac{\overline{y_i^2}-\overline{y_i}^2}{n}} = \pm 0.5 
$$
and the imaginary one 
$$
\Im(\hat{m})
=
\pm \hat{a} 
\sqrt{\overline{x_i^2}-\overline{x_i}^2}
= \mp 0.2 \ .
$$

\end{document}